\documentclass[12pt]{amsart}
\usepackage[all]{xy}
\usepackage{amssymb}
\usepackage{amsmath}
\usepackage{amsfonts}
\usepackage{amsthm}
\DeclareMathSizes{11.1}{10}{8}{6}

\newcommand{\simto}{\buildrel{\sim}\over{\rightarrow}}
\newcommand{\ra}{\rightarrow}
\renewcommand{\O}{{\mathcal O}}
\newcommand{\sub}{\subset}
\newcommand{\End}{\operatorname{End}}
\newcommand{\Z}{{\mathbb Z}}
\newcommand{\PP}{{\mathbb P}}
\newcommand{\Pic}{\operatorname{Pic}}
\newcommand{\Sing}{\operatorname{Sing}}
\newcommand{\Sym}{\operatorname{Sym}}
\newcommand{\wt}{\widetilde}

\newcommand{\lrar}[1]{\begin{picture}(50,10)(-25,-5)
\put(-25,0){\vector(1,0){50}}
\put(0,5){\makebox(0,0)[b]{\mbox{$#1$}}}
\end{picture}}

\newcommand{\ldar}[1]{\begin{picture}(10,50)(-5,-25)
\put(0,25){\vector(0,-1){50}}
\put(5,0){\mbox{$#1$}}
\end{picture}}

\title{Torelli theorem via Fourier-Mukai
transform.}
\author{A.~Beilinson, A.~Polishchuk}

\begin{document}

\vspace*{1.5cm}

\maketitle
\sloppy

 We show that the Fourier
transform on the Jacobian of a
curve interchanges ``$\delta$ functions" on the
curve and the theta divisor. The Torelli
theorem is an immediate
consequence.

\section{Statement of the theorem.}

\subsection{}\label{pre} We live over an
algebraically closed base field
$k$. Let $J$ be an
abelian variety equipped with a principal
polarization $\theta: J\simto
J^{\circ}=\Pic^0 (J)$, so we have the
corresponding Fourier transform
$\mathcal F$ on the derived category of
quasi-coherent sheaves
$D(J,\mathcal O )$ (see \cite{Muk}).

Let $\Theta$ be the theta divisor in $J$. Notice that
$\Theta$ is defined up to translation, and any
non-trivial translation does not preserve
$\Theta$. So we may consider $\Theta$ as a
canonically defined algebraic variety equipped
with a $J$-torsor of embeddings
$j:\Theta\hookrightarrow J$; we call these
$j$'s \emph{standard embeddings}. Denote by
$\Theta^{ns}$ the open subset of smooth points
of $\Theta$. For a standard embedding $j$ let
$j^{ns}:\Theta^{ns}\hookrightarrow J$ be its
restriction to $\Theta^{ns}$.

Our $\Theta$ carries a canonical involution
$x\mapsto x^{\nu}$; this is the unique
involution such that for any standard
embedding
$j$ the embedding $j^\nu :
x\mapsto -j(x^{\nu})$ is also standard.
For a line bundle $L$ on $\Theta$ or
$\Theta^{ns}$  set $ L^\nu :=\nu ^*
L$.

The pull-back $j^* F$ of an $\mathcal O_J$-module $F$
does not change if we translate
both $j$ and $F$ by the same element of $J$.
Thus the image of $j^{ns *}: \Pic (J) \to \Pic
(\Theta^{ns})$ is a canonically defined
subgroup of $\Pic (\Theta^{ns})$ (it does not
depend on $j$). Denote by $A(J)$ the
corresponding quotient group.

 Let $\mathcal T\subset
\Pic (\Theta^{ns})$ be the subset of line
bundles
$L$ such that

(i) $L\cdot L^\nu =\omega_{\Theta^{ns}}$

(ii) $A(J)$ is generated by the image of $L$.

\noindent
\emph{Remark.} Since the tangent bundle to $J$
is trivial, one has
$\omega_{\Theta^{ns}}=j^{ns*}\mathcal O_J
(j(\Theta ))$. Thus, if $\mathcal T$ is
non-empty then $\nu$ acts on $A(J)$ as $-1$.

\subsection{} From now on we assume that $J$ is
the Jacobian of a smooth projective curve $C$
of genus $g\ge 2$ equipped with the canonical
polarization. There is a standard embedding
$i:C\hookrightarrow J$ defined up to
translation; the standard embeddings $i$ form
a $J$-torsor isomorphic to $\Pic^{-1}(C)$.

\subsection{Theorem.} The set $\mathcal T$ is
non-empty.
For any $L\in\mathcal T$ and a standard
embedding $j:\Theta \hookrightarrow J$ the
 Fourier transform
$\mathcal F (j^{ns}_* L)$ equals to
$(\pm i)_* (M )[1-g]$ where
$i: C\hookrightarrow J$ is a standard
embedding, $M$ is a line bundle of degree $g-1$
on $C$.

\section{Proof of the theorem.}

\subsection{} Let us start with some
preliminaries.

Consider $J$ as the moduli space\footnote{We consider $J$ as a plain
variety ignoring the stack structure} of line bundles of
degree 0 on
$C$. A line bundle $E$ of degree $1-g$ yields a
standard embedding $j=j_E :\Theta
\hookrightarrow J$. Namely, there is a
canonical morphism $\sigma :
{\Sym}^{g-1}C\to\Theta$ such that $j_E \sigma$
sends a divisor $D$ to $E(D)$. Thus the
$J$-torsor of $j$'s equals
${\Pic}^{1-g}C$.
The above $\sigma$ is an isomorphism over
$\Theta^{ns}$; we denote by $\alpha$ the inverse
open embedding $\Theta^{ns}\hookrightarrow
{\Sym}^{g-1}C$.

For a group $T$ equipped with an
involution
$\nu$ we denote by $\wt{T}$ the
corresponding semi-direct product of $T$ and
$\Z /2\Z$. Our
$\mathcal T$ carries a canonical action of
$\wt{J}$ (here
$\nu$ acts on $J$ as $-1$). Namely, an element
$l\in J$ acts as tensor product
by the line bundle $j^{ns*}\theta (l)$ (notice
that, since
$\theta (l)$ is translation invariant, this
line bundle does not depend on
$j$ and is
$\nu$-anti-invariant), and $\Z /2\Z$ acts by
$\nu$.

\subsection{Proposition.}\label{p1} The
$\wt{J}$-action on $\mathcal T$ is transitive.

We prove \ref{p1} in \ref{p1p}. Notice that for
$g=2$ we  have $\Theta=\Theta^{ns}=C$, so here
the proposition is clear.

The map
$(L,j)\mapsto j^{ns}_*L$ commutes with the
action of $\wt{J}$; here $\wt{J}$ acts on
$\mathcal O$-modules on $J$ by twists by
degree $0$ line bundles and $-1$ symmetry.
The Fourier transform interchanges twists by
degree $0$ line bundles and translations and
commutes with the $-1$ symmetry. Since the set
of embeddings $\pm i:C\hookrightarrow J$ is a
$\wt{J}$-torsor,
\ref{p1} implies that it suffices to prove our
theorem for a single pair $(L,j)$.

Take any pair $(i,M)$ where $i:C
\hookrightarrow J$ is a standard embedding,
$M\in \Pic^{g-1}C$. The
theorem follows immediately from the
involutivity property of $\mathcal {F}$ and
the following fact:

\subsection{Proposition.}\label{p2} One has
$\mathcal {F} (i_* M)[1]=j_{*}^{ns}L$ where
$j:=j_{M^{ -1}}$ and
$L\in\mathcal T$.

We prove \ref{p2} in \ref{p2p}.

For every point $x\in C$ and every $d\ge 2$ consider the embedding
$$a_{x}^d:\Sym^{d-1}C\ra\Sym^dC:D\mapsto
D+x.$$ Let us denote by $R_{x}^d$ the image of
$a_{x}^d$.

\subsection{Lemma.}\label{picard} For $d\ge 2$
there is an exact sequence of abelian groups
$$0\longrightarrow \Pic(J)\stackrel{(j\sigma
)^*}\longrightarrow
\Pic(\Sym^{d}C)\stackrel{\deg}{\longrightarrow}\Z\longrightarrow
0$$ where the homomorphism $\deg$ is normalized
by the condition that
$\deg(\O(R_{x}^d))=1$ for any $x\in C$.
If $\PP\sub\Sym^d C$ is a complete linear
system of positive dimension then
$\deg(L)=\deg(L|_{\PP})$.

\emph{Proof of the lemma\footnote{An alternative proof can be found in
\cite{C}}.}
For $d$ sufficiently large the statement is true
since $\Sym^dC$ is a projective bundle over $J$.
Let $\pi_d:C^d\ra\Sym^dC$ be the canonical projection.
Then we have
$\pi_d^*\O(R_{x}^d
)\simeq\O_C(x)\boxtimes\ldots\boxtimes\O_C(x)$.
Therefore, $\O(R_{x}^d )$ is ample and
$H^i(\Sym^dC,\O(-nR_{x}^d ))=0$ for $n>0$,
provided that $i\le 1$, $d\ge 2$ or $i\le 2$,
$d\ge3$ (this follows from the symmetric Kunneth formula proved
in \cite{SGA4}, exp. XVII, 5.5.34).
It follows from the Lefschetz theorem
for Picard groups (see \cite{SGA2}, exp.11
(3.12), exp.12 (3.4)) that the map
$a_{x}^{d*}:\Pic(\Sym^dC)\ra\Pic(\Sym^{d-1}C)$
is an isomorphism for $d>3$ and is an
embedding for $d=3$. It remains to check that
$a_{x}^{3*}$ is surjective. Let us denote by
$K\subset\Pic(C\times C)$ the subgroup of line
bundles $L$ such that the restrictions
$L|_{x\times C}$ and $L|_{C\times x}$ are trivial. Then there is an
isomorphism $u:\End(J)\simto K:\phi\mapsto (i_x\times\phi i_x)^*\mathcal P$
where $i_x:C\ra J$ is the embedding corresponding to $x$,
$\mathcal P$ 
is the Poincar\'e line bundle on $J\times J$. Let $\Pic^+(C\times C)$
be the subgroup of line bundles stable under the involution $(x_1,x_2)\mapsto
(x_2,x_1)$, let $K^+=K\cap\Pic^+(C\times C)$. Then $u$ induces an isomorphism
of $\End^+(J)$ with $K^+$ where $\phi\in\End^+(J)$ if and only if
$\phi$ is self-dual. Let $r:\Pic(C\times C)\ra K$ be the homomorphism given by
$r(F)=F\otimes [(F^{-1}|_{C\times
x})\boxtimes(F^{-1}|_{x\times C})]$. Then
$r(\Pic^+(C\times C))=K^+$ and we have the
following commutative diagram

\begin{equation}
\setlength{\unitlength}{0.20mm}
\begin{array}{ccccc}
\Pic(J) & \lrar{\sigma_2^*} \Pic(\Sym^2C)
\lrar{\pi_2^*} & \Pic^{+}(C\times C)  \\
\ldar{s} &         & \ldar{r} \\
\End^+(J) & \setlength{\unitlength}{0.70mm} \lrar{u} & K^+
\end{array}
\end{equation}
where $\sigma:\Sym^2C\ra J$ maps $D$ to $\O(D-2x)$, $s:\Pic(J)\ra\End^+(J)$
is the standard homomorphism $L\mapsto\phi_L$ where $\phi_L(a)=
t_a^*L\otimes L^{-1}$. Since $s$ is surjective it follows that
the composition $r\circ\pi_2^*\circ\sigma_2^*$ is surjective.
Thus, in proving that some line bundle $L\in\Pic(\Sym^2C)$ comes
from $\Pic(\Sym^3C)$ we may assume that $\pi_2^*L$ belongs
to the kernel of $r$. In other words, $\pi_2^*L\simeq L_1\boxtimes L_1$ 
for some line bundle $L_1$ on $C$.
The $S_2$-action on $L_1\boxtimes L_1$ either coincides with the standard one
or differs from it by $-1$. In accordance with this dichotomy
we equip $\wt{L}=L_1\boxtimes L_1\boxtimes L_1$ either with the standard
$S_3$-action or with the standard action twisted by the sign character.
Then if we consider $\wt{L}$ as a line bundle on $\Sym^3 C$ we have
$a_{x}^{3*}\wt{L}=L$.

\subsection{Lemma.}\label{hyper} Assume that $C$ is hyperelliptic. Let
$Q\subset\Sym^{g-1}C$ be the complement to the image of $\alpha$.
Then $Q$ is an irreducible divisor and $\deg(Q)=-2$.

\emph{Proof of the lemma.}

Let $\tau:C\rightarrow C$ be the hyperelliptic involution.
For every $d>1$ let us denote by $Q_d\subset\Sym^dC$ the reduced
effective divisor consisting of $D$ such that $D$ contains a divisor
of the form $x+\tau x$.
Note that $Q_2\simeq\PP^1$ while
$Q_d$ is just the image of $Q_2\times\Sym^{d-2}C$ under
the natural map $\Sym^2C\times\Sym^{d-2}C\ra\Sym^dC$, so it is irreducible.
We have $Q=Q_{g-1}$.
It is easy to check that
\begin{equation}\label{res1}
a_{x}^{d*}\O(Q_d)\simeq\O(Q_{d-1}+R_{\tau(x)}^{d-1}).
\end{equation}
On the other hand, for any points $x,y\in C$ one has
\begin{equation}\label{res2}
a_{x}^{d*}\O(R_{y}^d )\simeq\O(R_{y}^{d-1}).
\end{equation}
Consider the embedding $a:\PP^1\simeq Q_2\hookrightarrow\Sym^dC$ given by
$D\mapsto D+D_0$ where $D_0=x_1+\ldots+x_{d-2}$ is a fixed effective divisor 
of degree $d-2$. Then by induction we derive from (\ref{res1}) and (\ref{res2})
that $a^*\O(R_{y}^d )\simeq\O_{\PP^1}(1)$ while
$a^*\O(Q_d)\simeq\O_{\PP^1}(Q_2\cdot Q_2+d-2)$, where $Q_2\cdot Q_2$ is the
self-intersection index of $Q_2$ in $\Sym^2C$. Since $Q_2\cdot Q_2=1-g$
we obtain $a^*\O(Q_d)\simeq\O_{\PP^1}(d-g-1)$.
In particular, $a^*\O(Q_{g-1})\simeq\O_{\PP^1}(-2)$ as required.

\subsection{Corollary.} Assume that $g\ge 3$.
Then
the map $j^{ns*}: \Pic (J) \to
\Pic (\Theta^{ns}
)$ is injective.
 The group $A(J)=\Pic
(\Theta^{ns})/j^{ns*}(\Pic (J))$ is isomorphic
to $\Z$ if $C$ is non-hyperelliptic, and to $\Z
/2\Z$ otherwise.

\emph{Proof of the corollary.}
(i) Assume that $C$ is
non-hyperelliptic.  Then by Martens theorem
the complement to the image of $\alpha$ has
codimension $>1$ (see \cite{ACGH}, IV 5.1)\footnote{The proof
given in \cite{ACGH} works in arbitrary
characteristic.}, so
$\alpha ^* : {\Pic}({\Sym}^{g-1}C)\simto {\Pic}(\Theta^{ns})$.
We are done by Lemma \ref{picard}.

(ii) Assume that $C$ is hyperelliptic. Then
Lemma
\ref{hyper} implies that $\Pic(\Theta^{ns})$
is isomorphic to $\Pic(\Sym^{g-1}C)/\Z[Q]$ where $\deg([Q])=2$,
so our statement follows easily from Lemma \ref{picard}.

\subsection{}\label{p1p}\emph{Proof of
Proposition
\ref{p1}.}

Choose $j$ to be
symmetric,
$j^{\nu}=j$;
then $j^{ns*}:{\Pic}J\to{\Pic}(\Theta^{ns})$ commutes
with the involution.
Take $L,L'\in\mathcal T$.
Since $L$ generates $A(J)$ and $L^{\nu}\equiv
L^{-1}\mod j^{ns*}\Pic(J)$ we have either
$L\equiv L'\mod j^{ns*}\Pic(J)$ or
$L^{\nu}\equiv L'\mod j^{ns*}\Pic(J)$.
Replacing $L$ by $L^{\nu}$ if necessary we
obtain that
$L^{-1}\cdot L'\simeq j^{ns*}\xi$ for some
$\xi\in\Pic(J)$.  Since $\xi^{\nu}=\xi^{-1}$
we deduce that $\xi\in\Pic^0(J)=J$ which
implies the proposition.

\emph{Remark.} We actually proved that if
$C$ is non-hyperelliptic then $\mathcal T$ is a
$\wt{J}$-torsor, while for hyperelliptic $C$,
it is a $J$-torsor (we did not prove that it
is non-empty as yet).

\subsection{}\label{p2p}\emph{Proof of
Proposition
\ref{p2}.} Our
$F:=\mathcal {F} (i_* M)[1]$ vanishes outside
of $j(\Theta)$ where $j=j_{M^{-1}}$. Since $F$ is the push-forward of
a (shifted) line bundle on $C\times J$ one may
represent it as the cone of a morphism
$f:V_1 \to V_0$ of vector bundles on $J$.
Therefore $f$ is injective (so $F=\operatorname{Coker} f$),
and for any closed subset
$Y\subset J$ of codimension
$>2$ one has $\mathcal {H}^1 _Y F=0$. 
Note also that $j(\Theta)$ is precisely the zero locus of $\det(f)$
(this follows from the well-known determinantal description of
$\Theta$, see e.g. \cite{ACGH}).
On the other hand, it is easy to see that $\det(f)$ annihilates 
$\operatorname{Coker} f$. Therefore, $F=j_*j^*F$. Since
the codimension in $J$ of the singular locus of $\Theta$ is $>2$
this shows that $F=j^{ns}_* L$ where $L:=j^{ns *}F$. By Riemann's theorem
on singularities of $\Theta$ the fibers of $L$ at all closed
points of $\Theta^{ns}$ are one-dimensional. 
Since $\Theta$ is reduced (see \cite{ACGH} IV.4.5) we see
that $L$ is a line bundle on $\Theta^{ns}$.

It remains to prove that
$L\in \mathcal T$. It is clear that the derived
pull-back
$\Phi :=Lj^{ns*} F$ is a complex
with
$H^0 \Phi =L$, $H^{-1}\Phi =L\otimes \mathcal
{O}_J (-j(\Theta ))$, other cohomology are 0.
On the other hand, computing
$\Phi$ using base change and the definition of
$F$ we see that for $x\in\Theta^{ns}$ the
corresponding fibers are
$H^0 \Phi_{j(x)} = H^1
(C,E_x \otimes M)$, $H^{-1}\Phi_{j(x)} = H^0
(C,E_x \otimes M)$ where 
$E_x :=i^* \theta(j(x))$ is
the line bundle on $C$ of degree 0 that
corresponds to $j(x)$.  Since 
$E_{x^\nu}\simeq E_x ^{-1}\otimes \omega_C \otimes M^{-2}$,  
the Serre duality yields
$(H^0\Phi_{j(x^{\nu})})^* = H^{-1}\Phi_{j(x)}$.
Therefore,
$L^{\nu *}=L\otimes \mathcal {O}_J
(-j(\Theta ))$. We see that
condition (i) from \ref{pre} is satisfied.

Let us check condition (ii). We may assume that
$i$ corresponds to a line bundle $\O_C (-x)$,
$x\in C$. Consider the universal divisor
$\mathcal D\subset C\times \Sym^{g-1}C$.
Then the pull-back of the Poincar\'e line bundle on $J\times J$
by the morphism $(i\times (j\sigma)):C\times\Sym^{g-1}C\ra J\times J$
is isomorphic to $p_1^*M^{-1}(\mathcal D-C\times
R_{x})$. It
follows that the line bundle $L^{-1}$ on
$\Theta^{ns}$ is
$\alpha^*p_{2*}(\O(\mathcal D))(-R_{x})$ where
$p_2:C\times\Sym^{g-1}C\ra\Sym^{g-1}C$ is the projection.
The canonical morphism $\O_{\Sym^{g-1}C}\ra p_{2*}(\O(\mathcal D))$
is an isomorphism over $\alpha(\Theta^{ns})$.
Therefore, $L^{-1}=\alpha^*\O(-R_{x})$ which generates
$A(J)$, so we are done.

\section{Concluding remarks.}

\subsection{} From
the Lefschetz theorem for Picard groups
(see \cite{SGA2}, exp.11 (3.12), exp.12 (3.4)) one can easily derive
that the restriction map $j^*:\Pic(J)\ra\Pic(\Theta)$ is an isomorphism
for an arbitrary principally polarized abelian variety $J$ of dimension
$g\ge4$. Furthermore, if the dimension of the singular locus
$\Sing\Theta$ is $<g-4$ then $\Theta$ is locally factorial
as follows from \cite{SGA2}, exp.11 (3.14). Hence, in this case
$\Pic(\Theta)=\Pic(\Theta^{ns})$ (since the notions of Cartier divisors
and Weil divisors on $\Theta$ coincide) and $A(J)=0$. Notice that if $J$ is a
Jacobian then the dimension of $\Sing\Theta$ is $\ge g-4$.
Moreover, Andreotti and Mayer proved in \cite{AM} that the closure of the
locus of Jacobians constitute an irreducible component of the locus
$N_{g-4}$ of principally polarized abelian varieties with $\dim\Sing\Theta\ge
g-4$. In \cite{B} Beauville established that in the case $g=4$ the locus
$N_{0}$ has two irreducible components. He also proved (assuming that
the characteristic is zero) that
a generic point of $N_{0}$ which is not contained in the closure of the locus
of Jacobians corresponds to an abelian variety $J$ with $\Sing\Theta$
consisting of one ordinary double point (see \cite{B}, 7.5). 
It follows that the corresponding
group $A(J)$ is either zero or
isomorphic to $\Z$ and the involution acts on $A(J)$ as
identity. Therefore, in this case either $A(J)=0$ or
the set $\mathcal T$ is empty.
The natural question is whether in higher dimensions one still has
either $A(J)=0$ or $\mathcal T=\emptyset$ 
for principally polarized abelian varieties which
are not in the closure of the locus of Jacobians.

\bigskip

{\bf Acknowledgments.} We are grateful to H\'el\`ene Esnault for pointing
out some errors in the first draft of the paper.  
The second author would like to thank J.~Harris and O.~Debarre
for helpful conversations
and A.~Collino for pointing out the reference \cite{C}.
A.~P. is partially supported by NSF grant.

\end{document}